\documentclass[12pt]{article}
\usepackage[top=2.25cm, bottom=2.25cm, left=2.25cm, right=2.25cm]{geometry}
\usepackage{amsmath,amssymb,amsthm,mathtools}
\usepackage{microtype}
\usepackage{enumitem}
\usepackage{booktabs}
\usepackage[hidelinks]{hyperref}
\usepackage[nameinlink,capitalize]{cleveref}
\usepackage{xcolor}

\newtheorem{theorem}{Theorem}[section]
\newtheorem{lemma}[theorem]{Lemma}

\newtheorem{corollary}[theorem]{Corollary}
\newtheorem{claim}[theorem]{Claim}
\theoremstyle{definition}

\newtheorem{conjecture}[theorem]{Conjecture}
\theoremstyle{remark}

\title{On several judicious partitioning problems of graphs and hypergraphs}
\author{Siwei Lin\footnote{Email: linsw0710@163.com},~~ Qinghou Zeng\footnote{Research partially supported by National Key R\&D Program of China (Grant No. 2023YFA1010202) and National Natural Science Foundation of China (Grant No. 12371342). Email: zengqh@fzu.edu.cn (Corresponding author)}\\
{\small Center for Discrete Mathematics, Fuzhou University, Fujian, 350003, China}}
\date{}

\begin{document}

\maketitle

\begin{abstract}
\textit{Judicious partition problems} ask for vertex partitions of graphs or hypergraphs that optimize several quantities simultaneously. In this paper, we investigate judicious partitions of graphs and hypergraphs in which each class meets many edges and establish the following results.
\begin{enumerate}
\item
Every hypergraph with nonnegative edge weights and all edges
of size at most $2$ has a bipartition in which each vertex class
meets edges of total weight at least
$(w_1-\Delta_1)/2+2w_2/3$,
where $w_i$ is the total weight of the edges of size $i$, and
$\Delta_1$ is the maximum weight of an edge of size $1$.

\item
Every hypergraph with all edges of size at most $3$ has
a bipartition in which each vertex class meets at least
$m_1/2+3m_2/4+7m_3/8+o(m_1+m_2+m_3)$ edges, where $m_i$ is the number of edges of size $i$.

\item
For any integer $k\ge4$, every graph with
$m\ge2(k-1)(k-2)$ edges has a partition into $k$ classes,
each meeting at least $2m/(2k-1)$ edges. Both the coefficient of $m$ and the edge threshold are best possible.

\item
For all sufficiently large $m$, every $4$-uniform
hypergraph with $m$ edges has a partition into four classes,
each meeting at least $4m/7$ edges.
\end{enumerate}

The first three statements are conjectures of Bollob\'as and Scott, and the last is a conjecture of Bollob\'as and Thomason on 4-uniform hypergraphs. All of the above results are best possible.

%The first two results answer a question of Bollob\'as and Scott on weighted hypergraphs with all edges of size at most $2$ and establish their conjectured asymptotic bound for hypergraphs with all edges of size at most $3$, respectively. The third verifies another conjecture of Bollob\'as and Scott on $k$-partitions of graphs, while the fourth confirms a conjecture of Bollob\'as and Thomason for sufficiently large 4-uniform hypergraphs.
\end{abstract}



\section{Introduction}

Classical graph or hypergraph partitioning problems usually ask for a partition that optimizes a single quantity.  For example, the well-known Max-Cut problem asks for a bipartition
$V_1,V_2$ of a graph $G$ that maximizes $e(V_1,V_2)$. Edwards \cite{Edwards73, Edwards75} proved the essentially best possible improvement $e(V_1,V_2)\ge m/2+(\sqrt{8m+1}-1)/8$.
Bollob\'as and Scott \cite{BS2002} extended this result to $k$-partitions by showing that every graph with $m$ edges admits a partition $V(G)=V_1\cup\cdots\cup V_k$ such that $e(V_1,\ldots,V_k) \ge \frac{k-1}{k}(m+(\sqrt{8m+1}-1)/4)-O(k)$.

While in practice, one often needs to find a partition of a given graph or hypergraph to optimize several quantities simultaneously. Such problems are called judicious partitioning problems by Bollob\'as and Scott. A basic example is the bottleneck bipartition problem, which was posed by Entringer \cite{SS1994}: given a graph $G$, find a bipartition $V_1,V_2$ of $V(G)$ that minimizes $\max\{e(V_1),e(V_2)\}$.  Thus one seeks a partition for which both $e(V_1)$ and $e(V_2)$ are small. Bollob\'as and Scott \cite{BS1999} showed that every graph with $m$ edges admits a bipartition $V_1\cup V_2$ such that 
$$\max\{e(V_1),e(V_2)\}\le\frac m4+\frac{\sqrt{8m+1}-1}{16}.$$
The bound is tight for the complete graph $K_{2n+1}$. 
In the same paper, the authors also extended the result for partitions into more than two parts. For more about judicious partition problems, we refer the readers to \cite{ABKS2003,BS1993,BS1997,BS2010,HWY2016,HX2027,JMYY2019,LSSZ2019,MY2010,MY2016,XYY2010,XY2011,XY2014}.

In this paper, we consider several types of judicious partitions problems for mixed $r$-graphs, i.e., hypergraphs with edges of size at most $r$. We first concern the bottleneck bipartition problem for mixed $r$-graphs.
Let $G$ be a mixed $r$-graph. We denote by $w_i(G)$ the total weight of edges of size $i$ for $1\le i\le r$ and denote by $\Delta_1$ the maximal weight of edges of size $1$.
 Bollob\'as and Scott \cite{BS2002} gave the following conjecture for weighted mixed 2-graphs.

\begin{conjecture}[Bollob\'as and Scott \cite{BS2002}]
Every weighted hypergraph \((G,w)\) with edges of size at most $2$ admits a bipartition such that each vertex class meets edges of total weight at least
\[
 \frac{w_1(G)-\Delta_1 }{2}+\frac{2w_2(G)}{3}.
\]
\end{conjecture}
Hou and Zeng \cite{HZ2017} established a weaker version of the conjecture. While, for weighted hypergraphs \(G\) with weight function \(w: E \to \mathbb{N}^+\), Haslegrave \cite{Haslegrave2012} confirmed the conjecture for the case \(\Delta_1 \le 1\). Here we state the following result for general weighted mixed $r$-graphs, which confirms the conjecture when $r=k=2$.

\begin{theorem}\label{thm:main k-part mixed graphs}
Let $r,k\ge2$ be integers. Let $G$ be a nonnegatively weighted hypergraph whose edges have size at most $r$. 
Then $V(G)$ has a $k$-partition $V_1\cup\ldots\cup V_k$ such that each class contains edges of total weight at most
\begin{equation*}
\sum_{j=2}^r\frac{w_j(G)}{j(k-1)+1}
 +\frac{w_1(G)+(k-1)\Delta_1}{k}.
\end{equation*}
\end{theorem}
\begin{corollary}
 Let $G$ be a nonnegatively weighted hypergraph whose edges have size at most $2$. Then $V(G)$ has a bipartition $V_1\cup V_2$ such that each class meets edges of total weight at least
\[
\frac{w_1(G)-\Delta_1}{2}+\frac{2w_2(G)}{3}.
\]
\end{corollary}

In the same paper, Bollob\'as and Scott \cite{BS2002} also suggested that the bound for nonweighted mixed $r$-graphs may be improved as follows when the edge number is sufficiently large.

\begin{conjecture}[Bollob\'as and Scott \cite{BS2002}]
Every hypergraph $H$ with $m_j$ edges of size $j$, $1\le j\le r$, has a bipartition $V(H)=V_1\cup V_2$ such that each class meets at least
\begin{align*}
\frac{m_1}{2}+\frac{3m_2}{4}+\frac{7m_3}{8}
 +\cdots+\left(1-\frac{1}{2^r}\right)m_r+o(m)
\end{align*}
edges, where $m=m_1+\cdots+m_r$.
\end{conjecture}
For $r=2$, Ma, Yen and Yu \cite{MYY2010} first confirmed this conjecture for the $k$-partitions version. Then Hou and Zeng \cite{HZ2017_2} gave the precise error term. For $r=3$, Tang, Yan and Zhang \cite{ZTY2015} established a weaker version of the conjecture. Our following result proves a general version.

\begin{theorem}\label{thm:main mixed 3-graphs}
Let $k\ge2$ be an integer, and let $p_1,\ldots,p_k>0$ satisfy $\sum_{i=1}^k p_i=1$.
Let $H$ be a hypergraph with edges of size at most $3$.  
Let $m_j$ denote the number of edges of size $j$ and $m=m_1+m_2+m_3$. Then there is a $k$-partition $V(H)=V_1\cup\cdots\cup V_k$ such that each class contains at most
\begin{align*}
 p_i m_1+p_i^2m_2+p_i^3m_3+O_k(m^{6/7})
\end{align*}
edges.
\end{theorem}
Setting $k=2$ and $p_1=p_2=1/2$, we confirm the conjecture of Bollob\'as and Scott for $r=3$.
\begin{corollary}
Let $H$ be a hypergraph with edges of size at most $3$.  
Let $m_j$ denote the number of edges of size $j$ and $m=m_1+m_2+m_3$.
Then there is a bipartition $V(H)=V_1\cup V_2$ such that each class meets at least
\begin{align*}
\frac{m_1}{2}+\frac{3m_2}{4}+\frac{7m_3}{8}+o(m)
\end{align*}
edges.
\end{corollary}

%We also consider another type of judicious partitioning problems, which asks for a partition in which each class meets many edges.
In \cite{BS2000}, Bollob\'as and Scott gave the following conjecture for graphs.

\begin{conjecture}[Bollob{\'a}s and Scott \cite{BS2000}]
For every integer \(k\ge2\), every graph with \(m\) edges has a partition into \(k\) sets, each of which meets at least
\[
 \frac{2m}{2k-1}
\]
edges.
\end{conjecture}

The complete graph \(K_{2k-1}\) shows that the coefficient \(2/(2k-1)\) cannot be improved. For $k=2$, the conjecture follows by considering a maximum bipartition of graphs. Zeng and Hou \cite{HZ2017} settled the case \(k=3\), apart from the two unavoidable small graphs \(K_2\) and \(K_{1,3}\). Here we confirm the conjecture for the edges number is at least $2(k-1)(k-2)$. It is worth mentioning that the star graph $K_{1,2(k-1)(k-2)-1}$ shows that the condition $m\ge2(k-1)(k-2)$ is necessary.

\begin{theorem}\label{thm:main k-part graphs}
Let \(k\ge4\) and let \(G\) be a simple graph with \(m\) edges. If $m\ge2(k-1)(k-2)$,
then \(V(G)\) has a partition \(V_1,\ldots,V_k\) such that each class meets at least
\[
 \left\lceil\frac{2m}{2k-1}\right\rceil.
\]
edges
\end{theorem}

For hypergraphs, Bollob\'as and Thomason (see \cite{BRT1993,BS2000}) gave a conjecture for uniform hypergraphs.
\begin{conjecture}[Bollob\'as and Thomason]
    For every integer \(r\ge2\), every $r$-uniform hypergraph with \(m\) edges has a partition into \(r\) sets, each of which meets at least
    \[
     \frac{rm}{2r-1}
    \]
    edges.
\end{conjecture}
For $r=3$, Bollob\'as and Scott \cite{BS2000} proved that every 3-uniform hypergraph with $m$ edges has a tripartition in which each class meets $\frac{5m-1}{9}$. Then Haslegrave \cite{Haslegrave2012} improved this bound to $0.6m$, which confirms the conjecture for $r=3$. Later, this bound was also improved asymptotically to $0.65m+o(m)$ by Ma and Yu \cite{MY2012}, and improved asymptotically to $\frac{19}{27}m+O(m^{6/7})$ by Spink and Tiba \cite{ST2020}. For general $r\ge3$, Bollob\'as, Reed and Thomason \cite{BRT1993} proved that every 3-uniform hypergraph with $m$ edges has a partition in which each class meets at least $(1-1/e )m/3\approx 0.21m$ edges, which was improved to $\frac{rm}{3r-4}$ by Haslegrave \cite{Haslegrave2014}.
Here we consider the conjecture for $r=4$ and prove the following conclusion.

\begin{theorem}\label{thm:main bipart 4-graphs}
Every simple $4$-uniform hypergraph $H$ with $m\ge1$ edges has a partition
$V(H)=V_1\cup V_2\cup V_3\cup V_4$ such that each class meets at least
\[
 \frac{37m}{64}-O(m^{8/9})
\]
edges.
\end{theorem}
Using this result, we confirm the conjecture of Bollob\'as and Thomason for $r=4$ and all sufficiently large $m$.
\begin{corollary}\label{coro:main bipart 4-graphs}
There is an absolute constant $m_0$ such that every simple $4$-uniform
hypergraph $H$ with $m\ge m_0$ edges has a partition
$V(H)=V_1\cup V_2\cup V_3\cup V_4$ such that each class meets at least
\[
\frac{4m}{7}
\]
edges.
\end{corollary}

\section{Preliminaries}
Let \( (H,w) \) be a weighted hypergraph with edges of size at most $r$. 
Let $i$ be an integer with $1\le i\le r$. An edge is an $i$-edge if its size is exactly $i$.
Let $E^{(i)}(H)$ be the set consisting of all $i$-edge of $H$.
For subsets $S,T$ of $V(H)$ and $i\ne 1$, $e^{(i)}_{(H,w)}(S,T)$ is the total weight of $i$-edges of \( H \) intersecting with both \( S \) and \( T \).
Let $e^{(i)}_{(H,w)}(S)$ be the total weight of $i$-edges of $H$ contained in $S$ and $w_i(H)=e^{(i)}_{(H,w)}(V(H))$.
For $i\ne1$, let \( d^{(i)}_{(H,w)}(S) \) be the total weight of $i$-edges of \( H \) meeting \( S \) (i.e., containing at least one vertex of \( S \)).
Further, we define $e_{(H,w)}(S,T)=\sum_{i\ne1}e^{(i)}_{(H,w)}(S,T)$, $e_{(H,w)}(S)=\sum_i e^{(i)}_{(H,w)}(S)$ and $d_{(H,w)}(S)=e^{(1)}_{(H,w)}(S)+\sum_{i\ne1} d^{(i)}_{(H,w)}(S)$.
If \( S = \{v\} \) and $i\ne 1$, then we write $w(v)$, \( e^{(i)}_{(H,w)}(v, T) \), $e_{(H,w)}(v, T)$, \( d^{(i)}_{(H,w)}(v) \) and $d_{(H,w)}(v)$ for $w(\{v\})$, \( e^{(i)}_{(H,w)}(\{v\}, T) \), $e_{(H,w)}(\{v\}, T)$, \( d^{(i)}_{(H,w)}(\{v\}) \) and $d_{(H,w)}(\{v\})$, respectively.
Let \( \Delta_i(H,w)=\max_{v}d^{(i)}_{(H,w)}(v) \) and $\Delta(H,w)=\max_{v}d_{(H,w)}(v)$.

For a nonweighted hypergraph $G$ with edges of size at most $r$, we similarly define $e^{(i)}_{G}(S,T)$, $e_{G}(S,T)$, $e^{(i)}_{G}(S)$, $e_{G}(S)$, $d^{(i)}_{G}(S)$, $d_{G}(S)$, \( e^{(i)}_{G}(v, T) \), $e_{G}(v, T)$, \( d^{(i)}_{G}(v) \), $d_{G}(v)$, $\Delta_i(G)$ and $\Delta(G)$.
Additionally, we write \( \overline{S} \) for \( V(G) \setminus S \), \([t]\) for \(\{1, \ldots, t\}\) and \( \binom{S}{j} \) for the set of all \( j \)-element subsets of \( S \). 

The subgraph of a hypergraph $G$ induced by a subset $U$ of $V(G)$, denoted by $G[U]$, is the hypergraph with vertex set $U$ whose edge set consists of all edges $e$ of $G$ with $e\subseteq U$.
Define the restriction of a weighted hypergraph $(G,w)$ to a subset $U$ of its vertices, denoted by $G|_U$, as the weighted hypergraph with vertex set $U$, in which the weight of an edge $e$ is the total weight of edges $f\in E(G)$ such that $e=f\cap U$.

\section{Judicious partitions of mixed hypergraphs}
In this section, we prove Theorem~\ref{thm:main k-part mixed graphs}.

\begin{proof}[\bf Proof of Theorem~\ref{thm:main k-part mixed graphs}.]
For $\ell \in\{2,\ldots,r\}$ and a $k$-partition
$\mathcal P=(V_1,\ldots,V_k)$, let
\[
 C_\ell(\mathcal P)
 \coloneqq
 w_\ell(G)-\sum_{i=1}^k e^{(\ell)}_{(G,w)}(V_i).
\]
Thus $C_\ell(\mathcal P)$ is the total weight of the $\ell$-edges that are not contained in a single part. Put $M\coloneqq w_1(G)+(k-1)\Delta_1(G,w)$.

We first assume that $M>0$. Define
\begin{equation}\label{eq:mixed-graph: function Phi}
 \Phi(\mathcal P)
 \coloneqq
 \sum_{\ell=2}^r
 \frac{C_\ell(\mathcal P)}{\ell(k-1)+1}
 +
 \frac1M
 \sum_{1\le a<b\le k}
 e^{(1)}_{(G,w)}(V_a)e^{(1)}_{(G,w)}(V_b).
\end{equation}
Choose a $k$-partition $\mathcal P=(V_1,\ldots,V_k)$ that is maximal for $\Phi$.
Fix $i\in[k]$, a vertex $v\in V_i$, and $j\in[k]\setminus\{i\}$.
Moving $v$ from $V_i$ to $V_j$ changes $C_\ell(\mathcal P)$ by
\[
\Delta C_\ell(v,\{i,j\})
\coloneqq
 \sum_{\substack{e\in E^{(\ell)}(G),\ v\in e\\
                   e\setminus\{v\}\subseteq V_i}}w(e)
 -
 \sum_{\substack{e\in E^{(\ell)}(G),\ v\in e\\
                   e\setminus\{v\}\subseteq V_j}}w(e).
\]
The change in the second term of \eqref{eq:mixed-graph: function Phi} is
\[
 \frac{w(v)}{M}
 \bigl(e^{(1)}_{(G,w)}(V_i)-e^{(1)}_{(G,w)}(V_j)-w(v)\bigr).
\]
The maximality of $\Phi$ therefore gives
\begin{align}\label{eq:mixed-single-move}
 -\sum_{\ell=2}^r\frac{1}{\ell(k-1)+1}
\Delta C_\ell(v,\{i,j\})\ge
 \frac{w(v)}{M}
 \bigl(e^{(1)}_{(G,w)}(V_i)-e^{(1)}_{(G,w)}(V_j)-w(v)\bigr).
\end{align}
For $\ell\in\{2,\ldots,r\}$, put
\[
 B_{\ell,i}
 \coloneqq
 \sum_{v\in V_i}\sum_{j\in[k]\setminus\{i\}}
 \sum_{\substack{e\in E^{(\ell)}(G),\ v\in e\\
                   e\setminus\{v\}\subseteq V_j}}w(e).
\]
Thus $B_{\ell,i}$ is the total weight of the $\ell$-edges having exactly one vertex in $V_i$ and all their remaining vertices in one other
part.
Sum \eqref{eq:mixed-single-move} first over
$j\in[k]\setminus\{i\}$ and then over $v\in V_i$. Since every
$\ell$-edge contained in $V_i$ is counted $\ell(k-1)$ times, we obtain
\begin{align}
 &\sum_{\ell=2}^r
 \frac{B_{\ell,i}-\ell(k-1)e^{(\ell)}_{(G,w)}(V_i)}
      {\ell(k-1)+1}
 \notag\\
 &\qquad\ge
 \frac1M
 \left(
 k\bigl(e^{(1)}_{(G,w)}(V_i)\bigr)^2
 -w_1(G)e^{(1)}_{(G,w)}(V_i)
 -(k-1)\sum_{v\in V_i}w(v)^2
 \right).
 \label{eq:mixed-local-master}
\end{align}
Since $0\le w(v)\le\Delta_1(G,w)$,
\[
 \sum_{v\in V_i}w(v)^2
 \le
 \Delta_1(G,w)e^{(1)}_{(G,w)}(V_i).
\]
Consequently,
\begin{equation}\label{eq:mixed-local-simplified}
 \sum_{\ell=2}^r
 \frac{B_{\ell,i}-\ell(k-1)e^{(\ell)}_{(G,w)}(V_i)}
      {\ell(k-1)+1}
 \ge
 \frac{e^{(1)}_{(G,w)}(V_i)
       \bigl(ke^{(1)}_{(G,w)}(V_i)-M\bigr)}{M}.
\end{equation}
For every $\ell$, the $\ell$-edges counted by $e^{(\ell)}_{(G,w)}(V_i)$ and those counted by $B_{\ell,i}$ are disjoint.
Hence
\[
 w_\ell(G)\ge e^{(\ell)}_{(G,w)}(V_i)+B_{\ell,i}.
\]
It follows that
\begin{align*}
 \sum_{\ell=2}^r\frac{w_\ell(G)}{\ell(k-1)+1}
 &\ge
 \sum_{\ell=2}^r
 \frac{e^{(\ell)}_{(G,w)}(V_i)+B_{\ell,i}}{\ell(k-1)+1}\\
 &=
 \sum_{\ell=2}^r e^{(\ell)}_{(G,w)}(V_i)
 +
 \sum_{\ell=2}^r
 \frac{B_{\ell,i}-\ell(k-1)e^{(\ell)}_{(G,w)}(V_i)}
      {\ell(k-1)+1}\\
 &\ge
 \sum_{\ell=2}^r e^{(\ell)}_{(G,w)}(V_i)
 +
 \frac{e^{(1)}_{(G,w)}(V_i)
       \bigl(ke^{(1)}_{(G,w)}(V_i)-M\bigr)}{M}.
\end{align*}
Therefore
\begin{align*}
 e_{(G,w)}(V_i)
 &\le
 \sum_{\ell=2}^r\frac{w_\ell(G)}{\ell(k-1)+1}
 +e^{(1)}_{(G,w)}(V_i)
 -\frac{e^{(1)}_{(G,w)}(V_i)
       \bigl(ke^{(1)}_{(G,w)}(V_i)-M\bigr)}{M}\\
 &=
 \sum_{\ell=2}^r\frac{w_\ell(G)}{\ell(k-1)+1}
 +2e^{(1)}_{(G,w)}(V_i)
 -\frac{k\bigl(e^{(1)}_{(G,w)}(V_i)\bigr)^2}{M}\\
 &=
 \sum_{\ell=2}^r\frac{w_\ell(G)}{\ell(k-1)+1}
 +\frac{M}{k}
 -\frac{k}{M}
 \left(e^{(1)}_{(G,w)}(V_i)-\frac{M}{k}\right)^2\\
 &\le
 \sum_{\ell=2}^r\frac{w_\ell(G)}{\ell(k-1)+1}
 +\frac{M}{k}.
\end{align*}

It remains to consider $M=0$. In this case
$w_1(G)=\Delta_1(G,w)=0$, and hence
$e^{(1)}_{(G,w)}(V_i)=0$ for every $i$. Choose a partition that is
locally maximal for
\[
 \sum_{\ell=2}^r
 \frac{C_\ell(\mathcal P)}{\ell(k-1)+1}.
\]
The same single-vertex variation argument, now without the quadratic
term, gives
\[
 \sum_{\ell=2}^r
 \frac{B_{\ell,i}-\ell(k-1)e^{(\ell)}_{(G,w)}(V_i)}
      {\ell(k-1)+1}
 \ge0.
\]
Therefore
\[
 e_{(G,w)}(V_i)
 =\sum_{\ell=2}^r e^{(\ell)}_{(G,w)}(V_i)
 \le
 \sum_{\ell=2}^r\frac{w_\ell(G)}{\ell(k-1)+1},
\]
which completes the proof.
\end{proof}

\section{Judicious partitions of mixed 3-graphs}

In this section, we prove Theorem \ref{thm:main mixed 3-graphs}. We first establish the convexity of the following implicit function.

\begin{lemma}\label{lem:root-convexity}
Let $a_0,T>0$ and $b_0,c,\alpha,\beta\ge0$.  Set
\begin{align*}
 a_t=a_0+\alpha t,\qquad b_t=b_0+\beta t
\end{align*}
for $t\ge0$.  For each $t$, let $r(t)>0$ be the unique solution of
\begin{align}\label{eq:cubic-root}
 a_tr+b_tr^2+cr^3=T.
\end{align}
Then $r$ is convex.
\end{lemma}

\begin{proof}[\bf Proof.]
For fixed $t$, the left-hand side of \eqref{eq:cubic-root}, viewed as a
function of $r\ge0$, is continuous, strictly increasing, vanishes at $r=0$,
and tends to infinity.  Thus the positive root exists and is unique.  Implicit
differentiation gives
\begin{align*}
 r'=-\frac{r(\alpha+\beta r)}{a_t+2b_tr+3cr^2}.
\end{align*}
A second differentiation gives
\begin{align*}
 r''=
 \frac{2r(\alpha+\beta r)}{(a_t+2b_tr+3cr^2)^3}
 \left[
   \alpha(a_t+b_tr)
   +\beta r(2a_t+3b_tr+3cr^2)
 \right].
\end{align*}
Every factor on the right-hand side is nonnegative, so $r''\ge0$.
\end{proof}

We shall also use Jensen's inequality in the following coordinatewise form.

\begin{lemma}\label{lem:iterated-jensen}
Let $n\ge1$, let $f:[0,1]^n\to\mathbb R$ be convex in each coordinate separately,
and let $X_1,\dots,X_n$ be independent Bernoulli random variables with
$\mathbb P(X_i=1)=\theta_i\in[0,1]$.  Then
\begin{align*}
 \mathbb E f(X_1,\dots,X_n)\ge f(\theta_1,\dots,\theta_n).
\end{align*}
\end{lemma}

\begin{proof}[\bf Proof.]
For $j=0,1,\dots,n$, set
\begin{align*}
 Y^{(j)}=(\theta_1,\dots,\theta_j,X_{j+1},\dots,X_n).
\end{align*}
Fix $j\in\{1,\dots,n\}$ and condition on $X_{j+1},\dots,X_n$.  Convexity in
the $j$th coordinate gives
\begin{align*}
 \mathbb E\!\left[f(Y^{(j-1)})\mid X_{j+1},\dots,X_n\right]
 \ge f(Y^{(j)}).
\end{align*}
Taking expectations, we obtain
\begin{align*}
 \mathbb E f(Y^{(j-1)})\ge \mathbb E f(Y^{(j)}).
\end{align*}
Iterating from $j=1$ to $j=n$ proves the result.
\end{proof}

We now prove an important inequality.  Let $D$ be a finite set.  For
$d\in D$, let $A_d,B_d\ge0$, and for each unordered pair
$\{d,d'\}\in\binom D2$, let $C_{dd'}\ge0$.  Write
\begin{align*}
 A=\sum_{d\in D}A_d,\qquad
 B=\sum_{d\in D}B_d,\qquad
 C=\sum_{\{d,d'\}\in\binom D2}C_{dd'}.
\end{align*}
For $S\subseteq D$, put
\begin{align*}
 A(S)=\sum_{d\in S}A_d,\qquad
 B(S)=\sum_{d\in S}B_d,\qquad
 C(S)=\sum_{\{d,d'\}\in\binom S2}C_{dd'}.
\end{align*}

\begin{lemma}\label{lem:multi-threshold}
Let $k\ge2$, let $\boldsymbol p=(p_1,\ldots,p_k)$ satisfy
$p_i>0$ and $\sum_{i=1}^k p_i=1$, and let
$\ell_1,\ell_2,\ell_3\ge0$.  Let the coefficients
$A_d,B_d,C_{dd'}$ and the functions $A(S),B(S),C(S)$ be as above.  For every $i\in[k]$, define
\begin{align*}
 T_i=p_i\ell_1+p_i^2\ell_2+p_i^3\ell_3
     +p_i^2A+p_i^3B+p_i^3C.
\end{align*}
Then there is an ordered partition
$D=S_1\cup\cdots\cup S_k$ and a probability
vector $\boldsymbol q=(q_1,\ldots,q_k)$ such that, for every $i\in[k]$,
\begin{align}\label{eq:multi-threshold-bound}
 \ell_1q_i+\ell_2q_i^2+\ell_3q_i^3
 +q_iA(S_i)+q_i^2B(S_i)+q_iC(S_i)
 \le T_i.
\end{align}
\end{lemma}

\begin{proof}[\bf Proof.]
If $D=\varnothing$, take $q_i=p_i$ for every $i$; then all inequalities hold
with equality.  Hence assume $D\neq\varnothing$.

Fix $\varepsilon>0$.  For $\mathbf x=(x_d)_{d\in D}\in[0,1]^D$, define
\begin{align*}
 a_\varepsilon(\mathbf x)
 =\ell_1+\varepsilon+\sum_{d\in D}A_dx_d
   +\sum_{\{d,d'\}\in\binom D2}C_{dd'}x_dx_{d'},
   \qquad
 b(\mathbf x)
 =\ell_2+\sum_{d\in D}B_dx_d.
\end{align*}
For each $i$, let $R_{i,\varepsilon}(\mathbf x)>0$ be the unique solution of
\begin{align}\label{eq:multi-root}
 a_\varepsilon(\mathbf x)R_{i,\varepsilon}(\mathbf x)+b(\mathbf x)R_{i,\varepsilon}(\mathbf x)^2
 +\ell_3R_{i,\varepsilon}(\mathbf x)^3=T_i+\varepsilon p_i.
\end{align}
The term $\varepsilon$ makes the coefficient of $R_{i,\varepsilon}$ positive throughout the
cube, so the left-hand side is strictly increasing in the root.

Fix all coordinates except $x_d$, and write $x_d=t$.  Then $a_\varepsilon(\mathbf x)$ is
affine in $t$ with nonnegative slope
\begin{align*}
 \alpha_d=A_d+\sum_{d'\ne d}C_{dd'}x_{d'}\ge0,
\end{align*}
while $b(\mathbf x)$ is affine in $t$ with nonnegative slope $\beta_d=B_d$.
Lemma~\ref{lem:root-convexity} shows that each $R_{i,\varepsilon}$ is convex in each
coordinate separately.

Color the elements of $D$ independently with colors $1,\ldots,k$, assigning
color $i$ with probability $p_i$, and let $S_i$ be the set of elements of color
$i$.  For fixed $i$, the coordinates of $\mathbf 1_{S_i}$ are independent
Bernoulli$(p_i)$ variables.  By Lemma~\ref{lem:iterated-jensen},
\begin{align}\label{eq:multi-jensen}
 \mathbb E R_{i,\varepsilon}(\mathbf 1_{S_i})
 \ge R_{i,\varepsilon}(p_i\mathbf 1_D).
\end{align}
At $\mathbf x=p_i\mathbf 1_D$,
\begin{align*}
 a_\varepsilon(\mathbf x)=\ell_1+\varepsilon+p_iA+p_i^2C,
 \qquad
 b(\mathbf x)=\ell_2+p_iB.
\end{align*}
Substituting $R=p_i$ into the left-hand side of \eqref{eq:multi-root} gives
$T_i+\varepsilon p_i$.  Therefore
\begin{align}\label{eq:multi-center-root}
 R_{i,\varepsilon}(p_i\mathbf 1_D)=p_i.
\end{align}
Summing \eqref{eq:multi-jensen} over $i$ and using
\eqref{eq:multi-center-root}, we get
\begin{align*}
 \mathbb E\sum_{i=1}^k R_{i,\varepsilon}(\mathbf 1_{S_i})
 \ge\sum_{i=1}^kp_i=1.
\end{align*}
Thus some deterministic coloring, denoted
$D=S_{1,\varepsilon}\cup\cdots\cup S_{k,\varepsilon}$, satisfies
\begin{align*}
 L_\varepsilon:=\sum_{i=1}^kR_{i,\varepsilon}(\mathbf 1_{S_{i,\varepsilon}})\ge1.
\end{align*}
For every $i\in[k]$, set
%\begin{align*}
 $q_{i,\varepsilon}
 =\frac{R_{i,\varepsilon}(\mathbf 1_{S_{i,\varepsilon}})}{L_\varepsilon}$.
%\end{align*}
Then $\sum_i q_{i,\varepsilon}=1$ and
$0\le q_{i,\varepsilon}\le R_{i,\varepsilon}(\mathbf 1_{S_{i,\varepsilon}})$.  Since the polynomial on
the left-hand side of \eqref{eq:multi-root} is increasing on $[0,\infty)$, for every $i\in[k]$,
\begin{align}
 (\ell_1+\varepsilon)q_{i,\varepsilon}+\ell_2q_{i,\varepsilon}^2+
 \ell_3q_{i,\varepsilon}^3
 +q_{i,\varepsilon}A(S_{i,\varepsilon})
 +q_{i,\varepsilon}^2B(S_{i,\varepsilon})
 +q_{i,\varepsilon}C(S_{i,\varepsilon})
 \le T_i+\varepsilon p_i.
 \label{eq:multi-regularized-bound}
\end{align}
Let $\varepsilon_n=1/n$.  There are only finitely many ordered partitions of $D$, so
after passing to a subsequence the sets $S_{i,\varepsilon_n}$ are constant.  By
compactness of the probability simplex, pass to a further subsequence for which
$(q_{1,\varepsilon_n},\ldots,q_{k,\varepsilon_n})$ converges.  Letting $n\to\infty$ in
\eqref{eq:multi-regularized-bound} proves \eqref{eq:multi-threshold-bound}.
\end{proof}

We use the following standard form of McDiarmid's bounded-difference
inequality \cite{McDiarmid1989}.

\begin{lemma}[McDiarmid's bounded-difference inequality]\label{lem:mcdiarmid}
Let $Z_1,\dots,Z_n$ be independent random variables, with $Z_i$ taking values
in a set $\Omega_i$, and let
$F:\Omega_1\times\cdots\times\Omega_n\to\mathbb R$.  Suppose that changing only
coordinate $i$ changes $F$ by at most $c_i$.  Then, for every $t>0$,
\begin{align*}
 \mathbb P(F(Z)-\mathbb E F(Z)\ge t)
 &\le \exp\!\left(-\frac{2t^2}{\sum_{i=1}^n c_i^2}\right),\\
 \mathbb P(F(Z)-\mathbb E F(Z)\le -t)
 &\le \exp\!\left(-\frac{2t^2}{\sum_{i=1}^n c_i^2}\right).
\end{align*}
If $\sum_i c_i^2=0$, then $F$ is constant and both conclusions are trivial.
\end{lemma}
Now we can prove \cref{thm:main mixed 3-graphs}.
\begin{proof}[\bf Proof of \cref{thm:main mixed 3-graphs}.]
The argument is trivial when $m=0$.  Let $\alpha=2/7$ and put
\begin{align*}
 U=\{v\in V(H):d_H(v)>m^{1-\alpha}\},
 \qquad
 W=V(H)\setminus U.
\end{align*}
Since
%\begin{align*}
 $\sum_{v\in V(H)}d_H(v)=m_1+2m_2+3m_3\le3m$,
%\end{align*}
we have
\begin{align}\label{eq:U-size}
 |U|\le3m^\alpha.
\end{align}
Let $\delta_j$ be the number of $j$-edges contained entirely in $U$, and put
$\delta=\delta_1+\delta_2+\delta_3$.  Since $H$ is simple,
\begin{align}\label{eq:delta-bound}
 e_H(U)
 &\le |U|+\binom{|U|}{2}+\binom{|U|}{3} \notag\\
 &=\left(\frac92+o(1)\right)m^{3\alpha}
  =\left(\frac92+o(1)\right)m^{6/7}.
\end{align}
Delete these edges temporarily, and call the resulting hypergraph $H'$.  Every
edge of $H'$ meets $W$.
For $j=1,2,3$, let $\ell_j$ be the number of $j$-edges contained in $W$.  For
$u\in U$, define
\begin{align*}
 A_u=\bigl|\{\{u,w\}\in E(H'):w\in W\}\bigr|,
 \qquad
 B_u=\bigl|\{\{u,w,w'\}\in E(H'):w,w'\in W\}\bigr|,
\end{align*}
and, for distinct $u,v\in U$, define
\begin{align*}
 C_{uv}=\bigl|\{\{u,v,w\}\in E(H'):w\in W\}\bigr|.
\end{align*}
Write
\begin{align*}
 A=\sum_{u\in U}A_u,
 \qquad
 B=\sum_{u\in U}B_u,
 \qquad
 C=\sum_{\{u,v\}\in\binom U2}C_{uv}.
\end{align*}
Thus $H'$ has $\ell_1$ one-edges, $\ell_2+A$ two-edges, and
$\ell_3+B+C$ three-edges.

Apply Lemma~\ref{lem:multi-threshold} with $D=U$ and
$\boldsymbol p=(p_1,\ldots,p_k)$.  We obtain an ordered partition $U=S_1\cup\cdots\cup S_k$ and a probability vector $(q_1,\ldots,q_k)$ such that, for every $i$,
\begin{align}\label{eq:expectation-target}
 \ell_1q_i+\ell_2q_i^2+\ell_3q_i^3
 +q_iA(S_i)+q_i^2B(S_i)+q_iC(S_i)
 \le T_i',
\end{align}
where
\begin{align*}
 T_i':=p_i\ell_1+p_i^2(\ell_2+A)+p_i^3(\ell_3+B+C).
\end{align*}
Place each vertex of $S_i$ in $V_i$.  Independently for each $w\in W$, place
$w$ in $V_i$ with probability $q_i$.  For $i\in[k]$, let $X_i=e_{H'}(V_i)$.
The six terms on the left-hand side of \eqref{eq:expectation-target} are
exactly the expected contributions to $X_i$ from the six edge types in $H'$.
Consequently, for $i\in[k]$,
\begin{align*}
 \mathbb E X_i\le T_i'.
\end{align*}
Moreover, changing the class of one vertex $w\in W$ changes any fixed $X_i$ by at most
$d_H(w)$.  Since $d_H(w)\le m^{1-\alpha}$ for $w\in W$,
\begin{align*}
 Q\coloneqq\sum_{w\in W}d_H(w)^2
 \le m^{1-\alpha}\sum_{w\in W}d_H(w)
 \le3m^{2-\alpha}.
\end{align*}
If $Q=0$, then $H'$ has no edges and the required simultaneous upper bounds are
immediate.  So we may suppose that $Q>0$.  By Lemma~\ref{lem:mcdiarmid}, for every $i$ and every
$z>0$,
\begin{align}\label{eq:tail-bound}
 \mathbb P(X_i>\mathbb E X_i+z)
 \le\exp\!\left(-\frac{2z^2}{Q}\right)
 \le\exp\!\left(-\frac{2z^2}{3m^{2-\alpha}}\right).
\end{align}
Choose
\begin{align}\label{eq:z-choice}
 z=\sqrt{\frac32\log(2k)}\;m^{1-\alpha/2}
  =\sqrt{\frac32\log(2k)}\;m^{6/7}.
\end{align}
Each probability in \eqref{eq:tail-bound} is at most $1/(2k)$.  A union bound
over the $k$ classes therefore shows that some realization satisfies
\begin{align}\label{eq:all-concentrated}
 X_i\le T_i'+z
\end{align}
 for every $i\in[k]$. The same conclusion is trivial when $Q=0$.

Restore all edges contained entirely in $U$.  For every $i$, the
internal-edge count increases by at most $e_H(U)$.  Moreover,
\begin{align*}
 T_i'
 =p_i(m_1-e^{(1)}_H(U))
   +p_i^2(m_2-e^{(2)}_H(U))
   +p_i^3(m_3-e^{(3)}_H(U))
   \le p_im_1+p_i^2m_2+p_i^3m_3.
\end{align*}
Using \eqref{eq:delta-bound}, \eqref{eq:z-choice}, and
\eqref{eq:all-concentrated}, we obtain, simultaneously for all $i$,
\begin{align*}
 e_H(V_i)
 \le T_i'+z+e_H(U)
 \le p_im_1+p_i^2m_2+p_i^3m_3
 +\left(\frac92+\sqrt{\frac32\log(2k)}+o(1)\right)m^{6/7}.
\end{align*}
After increasing $C$ if necessary to cover the finitely many small values of
$m$, this completes the proof.
\end{proof}

\section{Judicious partitions of graphs}
In this section, we prove \cref{thm:main k-part graphs}.
We divide the proof of \cref{thm:main k-part graphs} into two cases: a graph with large maximum degree and a graph with small maximum degree.
\subsection{Graphs with large maximum degree}

In this section we prove Theorem~\ref{thm:main k-part graphs} for graphs whose maximum degree is at least \(M\).
Write
\[
 (x)_+=\max\{x,0\}.
\]
For an integer-weighted hypergraph $(H,w)$ with edges of size at most $2$, define
\[
 g(H,w)=\sum_{v\in V(H)}\bigl(w(v)-1\bigr)_+,
 \qquad
 \Phi(H,w)=w_1(H)+w_2(H)-\frac12g(H,w).
\]

\begin{theorem}\label{thm:excess}
Let \(r\ge1\). Let \((H,w)\) be a weighted hypergraph with edges of
size at most \(2\), where \(w:E(H)\to\mathbb N^+\). If, for some positive integer \(M\), $\Delta(H,w)\le M-1$ and $\Phi(H,w)\ge rM$, then \(V(H)\) has an \(r\)-partition in which every class $X$ satisfies
\[
 d_{(H,w)}(X)\ge M.
\]
\end{theorem}

\begin{proof}[\bf Proof]
We induct on \(r\). If \(r=1\), then the whole vertex set meets all
edges, and
\[
 w_1(H)+w_2(H)\ge\Phi(H,w)\ge M.
\]
So we may suppose that \(r\ge2\). Since \(w_1(H)+w_2(H)\ge\Phi(H,w)\ge rM\), choose
\(S\subseteq V(H)\) satisfying
\[
 d_{(H,w)}(S)\ge M,
 \qquad
 |S|=\min\bigl\{|X|:X\subseteq V(H),\ d_{(H,w)}(X)\ge M\bigr\}.
\]
Write
\[
 d_{(H,w)}(S)=M+\epsilon-1,
 \qquad \epsilon\ge1,
 \qquad s=|S|.
\]
Since \(\Delta(H,w)\le M-1\), no singleton set meets total edge
weight at least \(M\), and hence \(s\ge2\).

For every \(v\in S\), the minimality of \(|S|\) gives $d_{(H,w)}(S\setminus\{v\})\le M-1$. Moreover,
\[
 d_{(H,w)}(S)-d_{(H,w)}(S\setminus\{v\})
 =d^{(1)}_{(H,w)}(v)+e_{(H,w)}(v,\overline S).
\]
Consequently,
\[
 w(v)+e_{(H,w)}(v,\overline S)\ge\epsilon.
\]
Summing over \(v\in S\), we obtain
\begin{align}\label{INEQ:private-weight-cardinality-choice}
 e^{(1)}_{(H,w)}(S)+e_{(H,w)}(S,\overline S)\ge s\epsilon.
\end{align}

Let \(H'=H|_{\overline S}\), and denote its weight function by $w'$. Define
\[
 N_1(S)=\bigl|\{v\in S:w(v)\ge1\}\bigr|,
 \quad
 N_2(\overline S)=
 \bigl|\{v\in\overline S:w(v)=0,
              \ e_{(H,w)}(v,S)\ge1\}\bigr|.
\]
For every \(v\in\overline S\), its weight in \(H'\) is $w'(v)=w(v)+e_{(H,w)}(v,S)$.
Hence
\[
 g(H',w')-g(H,w)
 =e_{(H,w)}(\overline S,S)-e^{(1)}_{(H,w)}(S)
   +N_1(S)-N_2(\overline S).
\]
Consequently,
\[
 \Phi(H,w)-\Phi(H',w')
 =e^{(2)}_{(H,w)}(S)+\frac12\bigl(
   e_{(H,w)}(S,\overline S)+e^{(1)}_{(H,w)}(S)
   +N_1(S)-N_2(\overline S)\bigr).
\]
Since
\[
 d_{(H,w)}(S)
 =e^{(1)}_{(H,w)}(S)+e^{(2)}_{(H,w)}(S)
  +e_{(H,w)}(S,\overline S),
\]
we obtain
\begin{align}\label{EQ:potential-change-under-restriction}
 d_{(H,w)}(S)-\bigl(\Phi(H,w)-\Phi(H',w')\bigr)=\frac12\bigl(
 e^{(1)}_{(H,w)}(S)+e_{(H,w)}(S,\overline S)
 -N_1(S)+N_2(\overline S)\bigr).
\end{align}
By \eqref{INEQ:private-weight-cardinality-choice},
\(N_1(S)\le s\), and \(N_2(\overline S)\ge0\), the right-hand side is
at least $\frac{s(\epsilon-1)}2\ge\epsilon-1$.
Thus
\[
 \Phi(H,w)-\Phi(H',w')
 \le d_{(H,w)}(S)-(\epsilon-1)=M,
\]
and therefore
\[
 \Phi(H',w')\ge(r-1)M.
\]
Moreover, restriction preserves the weighted degree of every surviving
vertex, so for each vertex $v\in V(H')$,
\[
 d_{(H',w')}(v)=d_{(H,w)}(v)\le M-1.
\]
The induction hypothesis partitions \(H'\) into \(r-1\) required
classes. Together with \(S\), these form the desired \(r\)-partition.
\end{proof}

\begin{theorem}\label{thm:high-degree}
Assume \(k\ge4\) and \(m\ge2(k-1)(k-2)\). If a graph \(G\) with \(m\)
edges has a vertex of degree at least
\(\left\lceil2m/(2k-1)\right\rceil\), then \(V(G)\) has a partition
\(V_1\cup\ldots\cup V_k\) such that, for each $i\in[k]$,
\[
 d_G(V_i)\ge\left\lceil\frac{2m}{2k-1}\right\rceil.
\]
\end{theorem}

\begin{proof}[\bf Proof]
Let
\[
 M=\left\lceil\frac{2m}{2k-1}\right\rceil,
 \qquad
 Q=\{v:d_G(v)\ge M\},
 \qquad
 q=|Q|.
\]
If \(q\ge k\), take any \(k\) vertices of \(Q\) as vertex classes and assign all remaining vertices arbitrarily. 
We may therefore assume \(1\le q\le k-1\), and we use every vertex of \(Q\) as one vertex class.

Let \(H=G|_{\overline Q}\), with weight function $w$. Then, for every $v\in \overline Q$,
%\[
 $w(v)=e_G(v,Q)\in\{0,1,\ldots,q\}$,
%\]
and every $2$-edge of $H$ has weight $1$. Thus for every vertex $v\in \overline Q$,
\[
 d_{(H,w)}(v)=d_G(v)\le M-1
\]
and
\begin{equation}\label{EQ:Total weight of H}
 w_1(H)+w_2(H)=m-e_G(Q).
\end{equation}
Write
\[
 K\coloneqq w_1(H)+w_2(H)-(k-q)M.
\]
By $m\ge(k-1)M$ and \(e_G(Q)\le\binom q2\),
\begin{align}\label{INEQ:K lower bound}
 K\ge(q-1)M-\frac{q(q-1)}2.
\end{align}

Let \(H_0=H\) and $w^0=w$. For each $i$, let $w^i$ denote the weight
function of $H_i$. We repeatedly choose an additional class. While fewer
than \(k-q\) additional classes have been constructed and
\[
 e^{(1)}_{(H,w)}(V(H_i))\ge M,
\]
choose \(C_{i+1}\) satisfying
\[
 C_{i+1}\subseteq
 \{v\in V(H_i):w(v)>0\},
 \qquad
 d_{(H_i,w^i)}(C_{i+1})\ge M,
\]
and
\[
 |C_{i+1}|=
 \min\Bigl\{|X|:
 X\subseteq\{v\in V(H_i):w(v)>0\},\
 d_{(H_i,w^i)}(X)\ge M\Bigr\}.
\]
Then put
\[
 H_{i+1}=H_i|_{V(H_i)\setminus C_{i+1}},
\]
and denote the weight function of $H_{i+1}$ by $w^{i+1}$.

We record the potential change in one step. Let \(C=C_{i+1}\), and define $a=e^{(1)}_{(H_i,w^i)}(C)$, $b=e^{(1)}_{(H,w)}(C)$, $c=e_{(H_i,w^i)}\left(C,V(H_i)\setminus C\right)$ and
\[
\begin{aligned}
 z=\left|\left\{v\in V(H_i)\setminus C:
 w^i(v)=0,\
 e_{(H_i,w^i)}(v,C)\ge1\right\}\right|.
\end{aligned}
\]
Write $d_{(H_i,w^i)}(C)=M+\epsilon-1$ where $\epsilon\ge1$.
For every \(v\in C\), the definition of \(C\) gives
\[
 d_{(H_i,w^i)}(C\setminus\{v\})\le M-1.
\]
Moreover, we have
\[
 d_{(H_i,w^i)}(C)-d_{(H_i,w^i)}(C\setminus\{v\})
 =w^i(v)
  +e_{(H_i,w^i)}\bigl(v,V(H_i)\setminus C\bigr).
\]
Thus
\[
 w^i(v)
 +e_{(H_i,w^i)}\bigl(v,V(H_i)\setminus C\bigr)\ge\epsilon,
\]
and summing over \(v\in C\) yields
\begin{align}\label{INEQ:current-private-weight-C}
 a+c\ge \epsilon|C|.
\end{align}
Furthermore, for each vertex $v\in C$, since $w(v)>0$, we have
\[
 w^i(v)\ge w(v)>0.
\]
The calculation leading to
\eqref{EQ:potential-change-under-restriction} therefore gives
\[
 \Phi(H_i,w^i)-\Phi(H_{i+1},w^{i+1})
 =e^{(2)}_{(H_i,w^i)}(C)+\frac12(a+c+|C|-z).
\]
Since
\[
 d_{(H_i,w^i)}(C)
 =a+e^{(2)}_{(H_i,w^i)}(C)+c=M+\epsilon-1,
\]
we obtain
\[
 M-\bigl(\Phi(H_i,w^i)-\Phi(H_{i+1},w^{i+1})\bigr)
 =\frac12(a+c-|C|+z)-(\epsilon-1).
\]
We claim that
\begin{align}\label{INEQ:Bound of difference of Phi for one step}
 \Phi(H_{i+1},w^{i+1})-\Phi(H_i,w^i)
 \ge\frac{b-|C|}{2}-M-(q-1).
\end{align}
Indeed, after subtracting the right-hand side of the equivalent
inequality, it remains to prove
\[
 \frac{a+c-b+z}{2}\ge\epsilon-q.
\]
If \(\epsilon\le q\), this follows from \(a\ge b\). If
\(\epsilon>q\), then \eqref{INEQ:current-private-weight-C} gives
\(a+c\ge \epsilon|C|\), while
\[
 b=e^{(1)}_{(H,w)}(C)\le q|C|.
\]
Hence, by \(z\ge0\) and \(|C|\ge2\),
\[
 a+c-b+z\ge (\epsilon-q)|C|\ge2(\epsilon-q).
\]
This proves \eqref{INEQ:Bound of difference of Phi for one step}.

If the procedure constructs all \(k-q\) additional classes, assign all
remaining vertices arbitrarily and finish. Otherwise it stops after
\(\ell\) steps, where
\[
 0\le\ell\le k-q-1,
 \qquad
 e^{(1)}_{(H,w)}(V(H_\ell))\le M-1.
\]
Let
\[
 F_0=\sum_{v\in \overline Q}\bigl(w(v)-1\bigr)_+,
 \qquad
 F_\ell=\sum_{v\in V(H_\ell)}
 \bigl(w(v)-1\bigr)_+.
\]
Every deleted vertex has positive base singleton weight, and the sets
\(C_1,\ldots,C_\ell\) are disjoint, so
\[
 F_0-F_\ell
 =\sum_{i=1}^\ell
 \bigl(e^{(1)}_{(H,w)}(C_i)-|C_i|\bigr).
\]
Summing \eqref{INEQ:Bound of difference of Phi for one step} and using
\(\Phi(H_0,w^0)=(k-q)M+K-F_0/2\), we obtain
\begin{align}\label{INEQ:Phi remaining high degree}
 \Phi(H_\ell,w^\ell)
 &\ge \Phi(H_0,w^0)+\frac12\sum_{i=1}^\ell
 \bigl(e^{(1)}_{(H,w)}(C_i)-|C_i|\bigr)
 -\ell M-\ell(q-1)\notag\\
 &=(k-q-\ell)M+K-\frac{F_\ell}{2}-\ell(q-1).
\end{align}
Since \(0\le w(v)\le q\), we have
%\[
 $\bigl(w(v)-1\bigr)_+
 \le\frac{q-1}{q}w(v)$,
%\]
and therefore
\begin{align}\label{INEQ:remaining base excess}
 F_\ell\le\frac{q-1}{q}(M-1).
\end{align}
Consequently, by $\ell\le k-q-1$,
\[
 \frac{F_\ell}{2}+\ell(q-1)
 \le\frac{(q-1)(M-1)}{2q}+(k-q-1)(q-1).
\]
Hence it is enough to prove
\begin{align}\label{INEQ:Deiffrence of K's bound}
 K\ge\frac{(q-1)(M-1)}{2q}+(k-q-1)(q-1).
\end{align}
By \eqref{INEQ:K lower bound}, the difference between its lower bound
for \(K\) and the right-hand side of
\eqref{INEQ:Deiffrence of K's bound} satisfies
\[
\begin{aligned}
 &(q-1)M-\frac{q(q-1)}2
 -\frac{(q-1)(M-1)}{2q}
 -(k-q-1)(q-1)\\
 &=\frac{1}{2q}(q-1)\bigl((2q-1)M+1-q(2k-q-2)\bigr)\\
 &\ge\frac{1}{2q}(q-1)\bigl((2q-1)(2k-4)+1-q(2k-q-2)\bigr)\\
 &=\frac{1}{2q}(q-1)^2(2k+q-5)\ge0,
\end{aligned}
\]
where
\[
 M=\left\lceil\frac{2m}{2k-1}\right\rceil
 \ge\left\lceil\frac{4(k-1)(k-2)}{2k-1}\right\rceil
 =2k-4
\]
since $m\ge2(k-1)(k-2)$. Thus
\eqref{INEQ:Deiffrence of K's bound} holds. Equations
\eqref{INEQ:Phi remaining high degree} and
\eqref{INEQ:remaining base excess} give
\[
 \Phi(H_\ell,w^\ell)\ge(k-q-\ell)M.
\]
Since $(H_\ell,w^\ell)$ is a restriction of $(H,w)$, foe each vertex $v\in V(H_\ell)$, $d_{(H_\ell,w^\ell)}(v)=d_{(H,w)}(v)$. So
\[
 \Delta(H_\ell,w^\ell)\le\Delta(H,w)\le M-1.
\]
Theorem~\ref{thm:excess} partitions \(H_\ell\) into the remaining
\(k-q-\ell\) classes. Together with the vertices of \(Q\) and the
sets \(C_1,\ldots,C_\ell\), this gives the required \(k\)-partition.
\end{proof}

\subsection{Graphs with small maximum degree}
In this section we prove Theorem~\ref{thm:main k-part graphs} for graphs whose vertices have degree at most
\(\left\lceil2m/(2k-1)\right\rceil-1\). We prove a more general weighted
statement.

\begin{theorem}\label{thm:capacity-q}
Let \(r\ge2\). Let \((H,w)\) be a weighted hypergraph with edges of
size at most \(2\), where \(w:E(H)\to\mathbb N^+\). If there is some
positive integer \(M\) such that $\Delta(H,w)\le M-1$ and
\[
 w_1(H)+2w_2(H)\ge(2r-1)M-2(r-1),
\]
then \(V(H)\) has an \(r\)-partition in which every class $X$ satisfies $d_{(H,w)}(X)\ge M$.
\end{theorem}

Since \(m\ge2(k-1)(k-2)\), together with $M=\left\lceil2m/(2k-1)\right\rceil$, gives
$2m\ge(2k-1)M-(2k-2)$,
the small-maximum-degree case follows from
Theorem~\ref{thm:capacity-q} by taking \(r=k\).

Next we first deal with the case $r=2$.

\begin{theorem}\label{thm:capacity-bip}
Let \((H,w)\) be a weighted hypergraph with edges of size at most \(2\), where \(w:E(H)\to\mathbb N^+\). If there is some positive integer \(M\) such that $\Delta(H,w)\le M-1$ and
\[
 w_1(H)+2w_2(H)\ge3M-2,
\]
then \(V(H)\) has a bipartition in which each class $X$ satisfies \(d_{(H,w)}(X)\ge M\).
\end{theorem}

We use the following singleton-pairing operation. If distinct vertices
\(u,v\) have positive singleton-edge weights, reduce the weights of
\(\{u\}\) and \(\{v\}\) by \(1\), and increase the weight of the
$2$-edge \(uv\) by \(1\), creating it if necessary. Denote the
transformed weighted hypergraph by \((H',w')\). Then we have the following conclusion.

\begin{claim}\label{clm:pairing}
Let \((H,w)\) and \((H',w')\) be as above. Then:
\begin{enumerate}[label=(\roman*)]
\item for every vertex $v\in V(H)$,
\[
 w_1(H)+2w_2(H)=w_1(H')+2w_2(H'),
 \qquad
 d_{(H,w)}(v)=d_{(H',w')}(v);
\]
\item for every \(X\subseteq V(H)\),
\[
 d_{(H,w)}(X)\ge d_{(H',w')}(X).
\]
\end{enumerate}
\end{claim}

\begin{proof}
A set containing neither \(u\) nor \(v\) is unaffected. A set containing
exactly one of them loses one unit of singleton-edge weight and gains one
unit from the edge \(uv\), so its meeting weight is unchanged. A set
containing both vertices loses two units of singleton-edge weight and gains
one unit from \(uv\), so its meeting weight decreases by one. This proves
(ii), while (i) follows immediately from the definition of the operation.
\end{proof}

To prove Theorem~\ref{thm:capacity-bip}, we need the following lemma.

\begin{lemma}\label{lem:partial-cover}
Let \((F,w)\) be a weighted hypergraph with edges of size at most
\(2\), where \(w:E(F)\to\mathbb N^+\). Let \(h,\ell\) be positive
integers. If
\begin{align}\label{INEQ:hypothesis of lemma partial cover}
 2w_2(F)\ge2h+w_1(F)+3\ell-2,
\end{align}
then there is a set \(U\subseteq V(F)\) such that
\begin{align}\label{INEQ:condition 1 of lemma partial cover}
 e^{(2)}_{(F,w)}(U)+e_{(F,w)}(U,\overline U)\ge\ell
\end{align}
and
\begin{align}\label{INEQ:condition 2 of lemma partial cover}
 e_{(F,w)}(U)
 =w_1(F)+w_2(F)-d_{(F,w)}(\overline U)
 \le w_2(F)-h-\ell.
\end{align}
\end{lemma}

\begin{proof}[\bf Proof]
Put $B=w_2(F)-h-\ell$, and define
\[
 \sigma=2w_2(F)-(2h+w_1(F)+3\ell-2)\ge0.
\]
Then we have
\begin{align}\label{EQ:relation of B and sigma}
 2B=w_1(F)+\ell-2+\sigma.
\end{align}
In particular, \(B\ge0\), since \(B\) is an integer and
\(2B\ge w_1(F)+\ell-2\ge-1\).

Add a new vertex \(x\) to \(F\). For every \(v\in V(F)\), give the
edge \(xv\) weight \(w(v)\). Denote the resulting
weighted graph by \((F',w')\). Choose a maximum-weight cut $(\{x\}\cup U,T)$
which, among all maximum-weight cuts, first minimizes $e_{(F,w)}(U)$ and then minimizes \(|U|\).

For \(v\in T\), moving \(v\) to \(\{x\}\cup U\) cannot increase the
cut weight. Hence
\[
 e_{(F,w)}(v,T\setminus\{v\})
 \le w(v)+e_{(F,w)}(v,U).
\]
Summing over \(v\in T\) gives
\[
 2e^{(2)}_{(F,w)}(T)
 \le w_1(F)-e^{(1)}_{(F,w)}(U)+e_{(F,w)}(U,T).
\]
Consequently,
\[
 3e^{(2)}_{(F,w)}(T)
 \le w_1(F)+e^{(2)}_{(F,w)}(U)+e_{(F,w)}(U,T)
       +e^{(2)}_{(F,w)}(T)
 =w_1(F)+w_2(F).
\]
Since \(h\ge1\), hypothesis
\eqref{INEQ:hypothesis of lemma partial cover} implies
\[
 2w_2(F)\ge w_1(F)+3\ell.
\]
Therefore
\[
 3e^{(2)}_{(F,w)}(T)
 \le w_1(F)+w_2(F)\le3w_2(F)-3\ell,
\]
so
\[
 e^{(2)}_{(F,w)}(T)\le w_2(F)-\ell.
\]
It follows that
\[
 e^{(2)}_{(F,w)}(U)+e_{(F,w)}(U,T)
 =w_2(F)-e^{(2)}_{(F,w)}(T)\ge\ell.
\]
Thus \(U\) satisfies \eqref{INEQ:condition 1 of lemma partial cover}.
If \(e_{(F,w)}(U)\le B\), it also satisfies
\eqref{INEQ:condition 2 of lemma partial cover}, and we are done.
Therefore we may suppose that $e_{(F,w)}(U)>B$ and set
\begin{align}\label{EQ:definition of xi}
 \xi:=e_{(F,w)}(U)-B\ge1.
\end{align}

Moving \(x\) to \(T\) cannot increase the cut weight. Hence $e^{(1)}_{(F,w)}(U)\le\frac{w_1(F)}2$.
Together with \eqref{EQ:relation of B and sigma},
and \eqref{EQ:definition of xi}, this yields
\begin{align}\label{INEQ:Precise bound of e(U)}
 2e^{(2)}_{(F,w)}(U)
 &=2e_{(F,w)}(U)-2e^{(1)}_{(F,w)}(U)\notag\\
 &\ge2B+2\xi-w_1(F)\notag\\
 &=\ell-2+\sigma+2\xi.
\end{align}
In particular, by $\sigma\ge0$ and $\xi\ge1$,
\[
 2e^{(2)}_{(F,w)}(U)\ge\ell.
\]

For every \(u\in U\), maximality of the cut gives
\[
 e_{(F,w)}(u,T)
 \ge w(u)+e_{(F,w)}(u,U\setminus\{u\}).
\]
If equality held, moving \(u\) to \(T\)
would preserve the cut weight and would either decrease \(e_{(F,w)}(U)\), or leave
\(e_{(F,w)}(U)\) unchanged while decreasing \(|U|\), contradicting the
lexicographic choice. So the inequality is strict. Since all weights are integers,
\begin{align}\label{INEQ:bound of e(u,T)}
 e_{(F,w)}(u,T)
 \ge w(u)
   +e_{(F,w)}(u,U\setminus\{u\})+1.
\end{align}
For \(Z\subseteq U\), summing
\eqref{INEQ:bound of e(u,T)} over \(u\in Z\) gives
\[
 e_{(F,w)}(Z,T)
 \ge e^{(1)}_{(F,w)}(Z)+2e^{(2)}_{(F,w)}(Z)
    +e_{(F,w)}(Z,U\setminus Z)+|Z|.
\]
Hence
\begin{align}\label{INEQ:bound of meeting edges of Z}
 e^{(2)}_{(F,w)}(Z)+e_{(F,w)}(Z,V(F)\setminus Z)\ge
 e^{(1)}_{(F,w)}(Z)+3e^{(2)}_{(F,w)}(Z)
 +2e_{(F,w)}(Z,U\setminus Z)+|Z|.
\end{align}
Define
\[
 L(Z)=e^{(1)}_{(F,w)}(Z)+3e^{(2)}_{(F,w)}(Z)
      +2e_{(F,w)}(Z,U\setminus Z)+|Z|.
\]

We now remove just enough of \(U\) to eliminate the excess cost
\(\xi\). For \(R\subseteq U\), put
\[
 \mu_U(R)=e_{(F,w)}(R)
          +e_{(F,w)}(R,U\setminus R).
\]
Since \(\mu_U(U)=e_{(F,w)}(U)\ge\xi\), choose \(R\subseteq U\) satisfying
\[
 \mu_U(R)\ge\xi,
 \qquad
 |R|=\min\bigl\{|X|:X\subseteq U,\ \mu_U(X)\ge\xi\bigr\}.
\]
Since \(\xi\ge1\), this implies that \(R\ne\varnothing\). Write
\[
 \mu_U(R)=\xi+\epsilon-1
\]
where $\epsilon\ge1$, and set \(Z=U\setminus R\). Removing \(R\) from \(U\) decreases
\(e_{(F,w)}(U)\) by exactly \(\mu_U(R)\).
Therefore
\begin{align}\label{EQ: Inside edge of Z}
 e_{(F,w)}(Z)
 =B-\epsilon+1.
\end{align}
Consequently,
\begin{align}\label{INEQ:Inside edge of Z}
 e_{(F,w)}(Z)\le B.
\end{align}
Thus \(Z\) satisfies \eqref{INEQ:condition 2 of lemma partial cover}.

It remains to verify \eqref{INEQ:condition 1 of lemma partial cover}.
For every \(v\in R\), the definition of \(R\) gives
\[
 \mu_U(R\setminus\{v\})\le\xi-1.
\]
Moreover,
\[
 \mu_U(R)-\mu_U(R\setminus\{v\})
 =w(v)+e_{(F,w)}(v,Z).
\]
Thus
\[
 w(v)+e_{(F,w)}(v,Z)\ge\epsilon.
\]
Summing over \(v\in R\), we obtain
\[
 e^{(1)}_{(F,w)}(R)+e_{(F,w)}(R,Z)\ge\epsilon|R|.
\]
Since $\mu_U(R)=e_{(F,w)}(R)+e_{(F,w)}(R,Z)$,
it follows that
\begin{align}\label{INEQ:lower bound of e(R)}
 e^{(2)}_{(F,w)}(R)=\mu_U(R)-\left(e^{(1)}_{(F,w)}(R)+e_{(F,w)}(R,Z)\right)\le\xi-1-\epsilon(|R|-1).
\end{align}
Directly from the definitions,
\[
 L(U)-L(Z)=\mu_U(R)+2e^{(2)}_{(F,w)}(R)+|R|.
\]
Using \(\mu_U(R)=\xi+\epsilon-1\) and
\eqref{INEQ:lower bound of e(R)}, we obtain
\begin{align*}
 L(U)-L(Z)
 &\le\xi+\epsilon-1
 +2\bigl(\xi-1-(|R|-1)\epsilon\bigr)+|R|\\
 &=3\xi+|R|+3\epsilon-3-2|R|\epsilon.
\end{align*}
On the other hand,
\[
 L(U)=e_{(F,w)}(U)+2e^{(2)}_{(F,w)}(U)+|U|.
\]
By \eqref{EQ:definition of xi} and
\eqref{INEQ:Precise bound of e(U)},
\[
 L(U)-\ell\ge B+3\xi+|U|+\sigma-2.
\]
Consequently,
\[
\begin{aligned}
 L(Z)-\ell
 &\ge B+3\xi+|U|+\sigma-2-\bigl(3\xi+|R|+3\epsilon-3-2|R|\epsilon\bigr)\\
 &=(B+1-\epsilon)+\sigma+(|U|-|R|)
   +2(|R|-1)\epsilon\\
 &\ge0.
\end{aligned}
\]
Indeed, \eqref{EQ: Inside edge of Z} gives
\[
 B+1-\epsilon=e_{(F,w)}(Z)\ge0,
\]
while \(\sigma\ge0\), \(|U|\ge|R|\), \(|R|\ge1\), and
\(\epsilon\ge1\). Hence \(L(Z)\ge\ell\), and
\eqref{INEQ:bound of meeting edges of Z} gives
\[
 e^{(2)}_{(F,w)}(Z)+e_{(F,w)}(Z,V(F)\setminus Z)\ge\ell.
\]
Together with \eqref{INEQ:Inside edge of Z}, this proves the lemma.
\end{proof}
Now we can prove \cref{thm:capacity-bip}.
\begin{proof}[\bf Proof of \cref{thm:capacity-bip}]
We repeatedly apply singleton pairing.

Let \(n=|V(H)|\), and list the vertices so that
\[
 w(v_1)\ge w(v_2)
 \ge\cdots\ge w(v_n).
\]
At every step, choose two vertices with the largest positive current
1-edge weights and apply the pairing operation. If $\Delta_1(H,w)\le w_1(H)/2$, then the analogous inequality remains valid after every step. If instead
\(w(v_1)>w_1(H)/2\), then the weight of
$\{v_1\}$ remains greater than half of the current total 1-edge
weight throughout the operation.

First suppose that $w(v_1)\le w_1(H)/2$. If \(w_1(H)\) is even, pair all singleton units across distinct vertices. 
If \(w_1(H)=2s+1\) is odd, then every 1-edge weight is at most \(s\), and at least three vertices have positive 1-edge weight. 
Here $s\ge1$. Reserve one singleton unit at three distinct positive vertices, chosen to include every vertex of 1-edge weight \(s\), and pair all remaining units across distinct vertices.
This is possible because the remaining total 1-edge weight is
\(2s-2\), and every remaining vertex has 1-edge weight at most
\(s-1\).

Let \((H',w')\) be the transformed weighted hypergraph. Then
\[
 w_1(H')\in\{0,3\},
\]
and \(\Delta_1(H',w')\le1\). By Claim~\ref{clm:pairing},
\[
 w_1(H')+2w_2(H')=w_1(H)+2w_2(H)\ge3M-2.
\]
Then Theorem~\ref{thm:main k-part mixed graphs} gives a bipartition of $H'$, in which each class meets at least
\[
 \frac{w_1(H')-\Delta_1(H',w')}{2}+\frac{2w_2(H')}{3}
 =\frac{w_1(H')+2w_2(H')}{3}
 \ge M-\frac23>M-1.
\]
Meeting weights are integers, so both classes meet at least \(M\).
Claim~\ref{clm:pairing} transfers the bipartition back to \((H,w)\).

Now suppose that $w(v_1)>w_1(H)/2$.
Pair every singleton unit outside \(v_1\) with a singleton unit at
\(v_1\). Let \((H',w')\) be the resulting weighted hypergraph. Only
\(v_1\) retains positive 1-edge weight, namely
\[
 h=w'(v_1)
   =2w(v_1)-w_1(H)>0.
\]
Let $F=H'|_{V(H')\setminus\{v_1\}}$, and denote its weight function by \(w_F\). By the definition of restriction,
\[
 w_1(F)=e_{(H',w')}(v_1,V(F)).
\]
Hence, by Claim~\ref{clm:pairing},
\[
 h+w_1(F)=d_{(H',w')}(v_1)=d_{(H,w)}(v_1)\le M-1.
\]
Define $\ell=M-h-w_1(F)$. Then \(\ell\ge1\). Moreover,
\[
 w_1(H')+2w_2(H')
 =h+2\bigl(w_1(F)+w_2(F)\bigr)\ge3M-2,
\]
so
\[
 2w_2(F)\ge2h+w_1(F)+3\ell-2.
\]
Apply Lemma~\ref{lem:partial-cover} to \((F,w_F)\), $h$ and $\ell$.
We obtain a set \(U\subseteq V(F)\) such that
\[
 e^{(2)}_{(F,w_F)}(U)+e_{(F,w_F)}(U,V(F)\setminus U)\ge\ell
\]
and
\[
 e_{(F,w_F)}(U)
 \le w_2(F)-h-\ell.
\]
Let \(V_1=V(F)\setminus U\). By the definition of $\ell$,
\[
\begin{aligned}
 d_{(H',w')}(V_1)
 =d_{(F,w_F)}(V_1)
 =w_1(F)+w_2(F)
   -e_{(F,w_F)}(U)
 \ge w_1(F)+h+\ell=M.
\end{aligned}
\]
Let \(V_2=U\cup\{v_1\}\). Since \(V_2\) meets the singleton edge at
\(v_1\), every edge incident with \(v_1\), and all edges of \(F\)
meeting \(U\),
\[
\begin{aligned}
 d_{(H',w')}(V_2)
 &\ge h+w_1(F)+e^{(2)}_{(F,w_F)}(U)
       +e_{(F,w_F)}(U,V_1)\\
 &\ge h+w_1(F)+\ell=M.
\end{aligned}
\]
Finally, Claim~\ref{clm:pairing} transfers this bipartition back to
\((H,w)\).
\end{proof}

We now prove Theorem~\ref{thm:capacity-q} by induction.

\begin{proof}[\bf Proof of Theorem~\ref{thm:capacity-q}]
We induct on \(r\). Since we deal with the case \(r=2\) in
Theorem~\ref{thm:capacity-bip}, we may assume that \(r\ge3\).

Suppose that the conclusion holds for $r-1$. The whole vertex set meets
at least \(M\) total edge weight. Indeed, if
\(w_1(H)+w_2(H)\le M-1\), then
\[
 w_1(H)+2w_2(H)\le2M-2,
\]
whereas
\[
 (2r-1)M-2(r-1)-(2M-2)
 =(2r-3)(M-1)+1>0,
\]
a contradiction.

Choose \(S\subseteq V(H)\) satisfying
\[
 d_{(H,w)}(S)\ge M,
 \qquad
 |S|=\min\bigl\{|X|:X\subseteq V(H),\ d_{(H,w)}(X)\ge M\bigr\}.
\]
Write
\[
 d_{(H,w)}(S)=M+\epsilon-1,
 \qquad \epsilon\ge1.
\]
For every \(v\in S\), the definition of $S$ gives
\[
 d_{(H,w)}(S\setminus\{v\})\le M-1,
\]
and
\[
 d_{(H,w)}(S)-d_{(H,w)}(S\setminus\{v\})
 =w(v)+e_{(H,w)}(v,\overline S).
\]
Consequently,
\[
 w(v)+e_{(H,w)}(v,\overline S)\ge\epsilon.
\]
Summing over \(v\in S\), we obtain
\[
 e^{(1)}_{(H,w)}(S)+e_{(H,w)}(S,\overline S)
 \ge\epsilon|S|.
\]
Since \(d_{(H,w)}(v)\le M-1\) for every vertex, \(|S|\ge2\), and hence
\[
\begin{aligned}
 e^{(2)}_{(H,w)}(S)
 =d_{(H,w)}(S)-e^{(1)}_{(H,w)}(S)
   -e_{(H,w)}(S,\overline S)
 \le M-\epsilon(|S|-1)-1
 \le M-\epsilon-1.
\end{aligned}
\]

Let \(H'=H|_{V(H)\setminus S}\), with weight function \(w'\). Then
\begin{align*}
 w_1(H')=w_1(H)-e^{(1)}_{(H,w)}(S)
                 +e_{(H,w)}(S,\overline S),
                 \quad
 w_2(H')=w_2(H)-e^{(2)}_{(H,w)}(S)
                 -e_{(H,w)}(S,\overline S).
\end{align*}
Therefore
\[
\begin{aligned}
 w_1(H')+2w_2(H')
 &=w_1(H)+2w_2(H)-d_{(H,w)}(S)-e^{(2)}_{(H,w)}(S)\\
 &\ge((2r-1)M-2(r-1))
   -(M+\epsilon-1)-(M-\epsilon-1)\\
 &=(2r-3)M-2(r-2).
\end{aligned}
\]
Restriction preserves every surviving weighted degree. Hence the
induction hypothesis gives an \((r-1)\)-partition of \((H',w')\) in
which every class meets at least \(M\). Together with \(S\), this is the
desired \(r\)-partition of \((H,w)\).
\end{proof}

\subsection{Proof of Theorem~\ref{thm:main k-part graphs}}

\begin{proof}[\bf Proof of Theorem~\ref{thm:main k-part graphs}]
Put \(M=\left\lceil2m/(2k-1)\right\rceil\). If \(G\) has a vertex of
degree at least \(M\), apply Theorem~\ref{thm:high-degree}. Otherwise
\(d_G(v)\le M-1\) for every vertex. Regard \(G\) as a weighted
hypergraph with no singleton edges and with every $2$-edge of weight
\(1\). Then
\[
 w_1(G)+2w_2(G)=2m\ge(2k-1)M-2(k-1)
\]
since $m\ge2(k-1)(k-2)$ and
$M=\left\lceil2m/(2k-1)\right\rceil$. Theorem~\ref{thm:capacity-q},
with \(r=k\), gives the required \(k\)-partition.
\end{proof}

\section{Judicious partitions of 4-uniform graphs}
In this section, we prove Theorem~\ref{thm:main bipart 4-graphs}. We first prove the following key lemma.

\begin{lemma}\label{lem:probabilities}
    Fix a nonnegative integer $m_0$. Let $F$ be a multihypergraph with $m_i$ edges of size at most $i$ for $i=1,2,3$. 
    There are a $4$-partition $V(F)=V_1\cup V_2\cup V_3\cup V_4$ and four reals $p_1,p_2,p_3,p_4\ge0$ with $p_1+p_2+p_3+p_4=1$ such that, for $i=1,2,3,4$,
    \begin{align*}
        (1-p_i)^4m_0+\sum_{e\cap V_i=\varnothing}(1-p_i)^{4-|e|}\le \frac{27}{64}(m_0+m_1+m_2+m_3).
    \end{align*}
\end{lemma}
\begin{proof}[\bf Proof.]
    Since it is trivial for $\sum_{i=0}^3m_i=0$, we may suppose that $\sum_{i=0}^3m_i>0$. First give the vertices of $V(F)$ independent uniform colors from $\{1,2,3,4\}$. 
    For $i=1,2,3,4$, let $U_i$ be the set of vertices with color $i$, and set
    \[
        s_i=\sum_{\substack{e\in E(F)\\e\cap U_i=\varnothing}}
          \left(\frac43\right)^{|e|}.
    \]
    Since an edge avoids color $i$ with probability $(3/4)^{|e|}$, we have $\mathbb E s_i=m_1+m_2+m_3$ and $\mathbb E\sum_i s_i=4(m_1+m_2+m_3)$.
    Fix a coloring $V_1\cup V_2\cup V_3\cup V_4$ of $V(F)$ for which $\sum_i s_i\le4(m_1+m_2+m_3)$, and put
    \[
        q_i=\frac{s_i+m_0}{4(m_0+m_1+m_2+m_3)},\qquad
        p_i=q_i+\frac{1-\sum_{j=1}^4q_j}{4}.
    \]
    Then $0\le q_i\le p_i\le1$ and $\sum_ip_i=1$.

    Let $Z_i=(1-p_i)^4m_0+\sum_{e\cap U_i=\varnothing}(1-p_i)^{4-|e|}$. Then
    \begin{align*}
        Z_i&=\left(\frac43\right)^{0}\left(\frac34\right)^{0}(1-p_i)^4m_0+\sum_{\substack{e\in E(F)\\e\cap U_i=\varnothing}}\left(\frac43\right)^{|e|}\left(\frac34\right)^{|e|}(1-p_i)^{4-|e|}\\
        &\le(s_i+m_0)\max_{0\le j\le3}\left\{\left(\frac34\right)^{j}(1-p_i)^{4-j}\right\}\\
        &\le(s_i+m_0)\max_{0\le j\le3}\left\{\left(\frac34\right)^{j}(1-q_i)^{4-j}\right\}\\
        &=4(m_0+m_1+m_2+m_3)\max\left\{q_i(1-q_i)^{4},\left(\frac34\right)^{3}q_i(1-q_i)\right\}.
    \end{align*}
    Note that for $0\le x\le1$, we have $x(1-x)^4\le\frac14(\frac45)^5$ and $\left(\frac34\right)^{3}x(1-x)\le\frac14\left(\frac34\right)^{3}$. 
    Thus we have
    \begin{align*}
        Z_i\le \frac{27}{64}(m_0+m_1+m_2+m_3).
    \end{align*}
    This completes the proof.
\end{proof}

Now we can prove Theorem~\ref{thm:main bipart 4-graphs}.

\begin{proof}[\bf Proof of Theorem~\ref{thm:main bipart 4-graphs}]
Let $D=2m^{7/9}$ and put
\[
A=\{v\in V(H):d_H(v)\ge D\},
\qquad U=V(H)\setminus A.
\]
Since $\sum_{v\in V(H)}d_H(v)=4m$, we have $|A|\le\frac{4m}{D}=2m^{2/9}$.
Consequently,
\[
e_H(A)\le\binom{|A|}{4}
\le\frac23m^{8/9}.
\]
Let $H'$ be obtained from $H$ by deleting all edges contained entirely in $A$.
Let $F$ be the restriction of $H'$ to $A$, that is, the multihypergraph with vertex set $A$, in which an edge $f$ of $F$ has the multiplicity $|\{e\in E(H'):e\cap A=f\}|$.
For $j=1,2,3$, let $m_j$ be the number of edges of $F$ of size $j$, counting multiplicities, and put $m_0=e_H(U)$. Then $M:=m_0+m_1+m_2+m_3=m-e_H(A)>0$.

By Lemma~\ref{lem:probabilities}, there exist a partition $A=A_1\cup A_2\cup A_3\cup A_4$ and nonnegative reals $p_1,p_2,p_3,p_4$ with $\sum_i p_i=1$ such that, for every $i\in\{1,2,3,4\}$,
\[
(1-p_i)^4m_0+
\sum_{\substack{f\in E(F)\\
f\cap A_i=\varnothing}}
(1-p_i)^{4-|f|}
\le\frac{27M}{64}.
\]
Independently assign each vertex of $U$ color $i$ with probability $p_i$.
Let $U_i$ be the set of vertices receiving color $i$, and put $V_i=A_i\cup U_i$.
For each $i$, let $X_i$ be the number of edges of $H'$ avoiding $V_i$.
Then
\[
\begin{aligned}
\mathbb E X_i
=(1-p_i)^4m_0+
\sum_{\substack{f\in E(F)\\
f\cap A_i=\varnothing}}
(1-p_i)^{4-|f|}\le\frac{27M}{64}.
\end{aligned}
\]
Changing the color of a single vertex $u\in U$ changes $X_i$ by at most $d_H(u)$. Set $t=\sqrt{\frac{\log8}{2}\sum_{u\in U}d_H(u)^2}$.
By Lemma~\ref{lem:mcdiarmid} (McDiarmid's inequality),
\[
\mathbb P(X_i\ge\mathbb EX_i+t)
\le\frac18.
\]
A union bound shows that there is an assignment for which, simultaneously for $i=1,2,3,4$, $X_i\le\frac{27M}{64}+t$.
Since every vertex of $U$ has degree less than $D$, we have
\[
t\le\sqrt{2(\log8)Dm}
=2\sqrt{\log8}\,m^{8/9}.
\]
Therefore, for $V_1\cup V_2\cup V_3\cup V_4$,
\[
\begin{aligned}
d_H(V_i)
\ge d_{H'}(V_i)
=M-X_i
\ge\frac{37m}{64}
-\left(\frac{37}{96}
+2\sqrt{\log8}\right)m^{8/9}.
\end{aligned}
\]
This holds for all $i\in\{1,2,3,4\}$,
completing the proof.
\end{proof}

\section*{Declaration on the Use of AI}

The authors used ChatGPT 5.6 to assist in discussing proof strategies, checking proofs, and improving exposition.


\begin{thebibliography}{99}
\bibitem{ABKS2003}
N.~Alon, B.~Bollob\'as, M.~Krivelevich and B.~Sudakov,
\emph{Maximum cuts and judicious partitions in graphs without short cycles},
J. Combin. Theory Ser. B \textbf{88} (2003) 329--346.

\bibitem{BRT1993}
B. Bollob\'as, B.~A. Reed and A. Thomason,
\textit{An extremal function for the achromatic number},
in {\it Graph structure theory (Seattle, WA, 1991)}, 161--165, Contemp. Math., 147, Amer. Math. Soc., Providence, RI, 1993.

\bibitem{BS1993}
B.~Bollob\'as and A.~D.~Scott,
\emph{Judicious partitions of graphs},
Period. Math. Hungar. \textbf{26} (1993) 125--137.

\bibitem{BS1997}
B.~Bollob\'as and A.~D.~Scott,
\emph{Judicious partitions of hypergraphs},
J. Combin. Theory Ser. A \textbf{78} (1997) 15--31.

\bibitem{BS1999}
B.~Bollob\'as and A.~D.~Scott,
\emph{Exact bounds for judicious partitions of graphs},
Combinatorica \textbf{19} (1999) 473--486.

\bibitem{BS2000}
B.~Bollob\'as and A.~D.~Scott,
\emph{Judicious partitions of 3-uniform hypergraphs},
European J. Combin. \textbf{21} (2000) 289--300.

\bibitem{BS2002}
B.~Bollob\'as and A.~D.~Scott,
\emph{Problems and results on judicious partitions},
Random Struct. Alg. \textbf{21} (2002) 414--430.

\bibitem{BS2010}
B.~Bollob\'as and A.~D.~Scott,
\emph{Max k-cut and judicious k-partitions},
Discrete Math. \textbf{310} (2010) 2126--2139.

\bibitem{Edwards73}
C.~S.~Edwards, 
\emph{Some extremal properties of bipartite subgraphs}, 
Canad. J. Math. \textbf{25} (1973) 475--485.

\bibitem{Edwards75}
C.~S.~Edwards,
\emph{An improved lower bound for the number of edges in a largest bipartite subgraph},
in \emph{Recent Advances in Graph Theory (Proc. Second Czechoslovak Sympos., Prague, 1974)},
Academia, Prague, (1975) 167--181.

\bibitem{FH2017}
G. Fan and J. Hou,
\textit{Bounds for pairs in judicious partitioning of graphs},
Random Struct. Alg. {\bf 50} (2017) 59--70.

\bibitem{Haslegrave2012}
J.~Haslegrave,
\emph{The Bollob\'as--Thomason conjecture for 3-uniform hypergraphs},
Combinatorica \textbf{32} (2012) 451--471.

\bibitem{Haslegrave2014}
J. Haslegrave,
\textit{Judicious partitions of uniform hypergraphs},
Combinatorica {\bf 34} (2014) 561--572.

\bibitem{HWY2016}
J. Hou, S. Wu and G. Yan,
\textit{On judicious partitions of uniform hypergraphs},
J. Combin. Theory Ser. A {\bf 141} (2016) 16--32.

\bibitem{HZ2017}
J.~Hou and Q.~Zeng,
\emph{Bounds for the number of meeting edges in graph partitioning},
Czechoslovak Math. J. \textbf{67} (2017) 741--752.

\bibitem{HZ2017_2}
J.~Hou and Q.~Zeng,
\emph{Judicious partitioning of hypergraphs with edges of size at most 2},
Combin. Probab. Comput. \textbf{26} (2017) 267--284.

\bibitem{HX2027}
Q.~Hu and B.~Xu,
\textit{A note on judicious partition of hypergraphs},
Discrete Math. \textbf{350} (2027) 115349.

\bibitem{JMYY2019}
Y.~Ji, J.~Ma, J.~Yan and X.~Yu, 
\textit{On problems about judicious bipartitions of graphs},
J. Combin. Theory Ser. B \textbf{139} (2019) 230--250.

\bibitem{LSSZ2019}
D. Lokshtanov, S. Saurabh, R. Sharma and M. Zehavi,
\textit{Balanced judicious bipartition is fixed-parameter tractable},
SIAM J. Discrete Math. {\bf 33} (2019) 1878--1911.

\bibitem{MYY2010}
J.~Ma, P.~Yen and X.~Yu,
\emph{On several partitioning problems of Bollob\'as and Scott},
J. Combin. Theory Ser. B \textbf{100} (2010) 631--649.

\bibitem{MY2010}
J. Ma and X. Yu,
\textit{Bounds for pairs in partitions of graphs},
Discrete Math. {\bf 310} (2010) 2069--2081.

\bibitem{MY2012}
J. Ma and X. Yu,
\textit{Partitioning 3-uniform hypergraphs},
J. Combin. Theory Ser. B {\bf 102} (2012) 212--232.

\bibitem{MY2016}
J. Ma and X. Yu,
\textit{On judicious bipartitions of graphs},
Combinatorica {\bf 36} (2016) 537--556.

\bibitem{McDiarmid1989}
C.~McDiarmid,
\newblock On the method of bounded differences,
\newblock in J.~Siemons (ed.), \emph{Surveys in Combinatorics, 1989},
London Math. Soc. Lecture Note Ser. 141, Cambridge Univ. Press, 1989,
pp.~148--188.

%\bibitem{Scott2005}
%A. Scott,
%\emph{Judicious partitions and related problems},
%in: Surveys in Combinatorics 2005,
%London Math. Soc. Lecture Note Ser. 327,
%Cambridge Univ. Press, Cambridge, 2005, 95--117.

\bibitem{SS1994}
F.~Shahrokhi and L.~A.~Sz\'ekely,
\emph{The complexity of the bottleneck graph bipartition problem},
J. Combin. Math. Combin. Comput. \textbf{15} (1994) 221--226.

\bibitem{ST2020}
H. Spink and M. Tiba,
\textit{Judiciously 3-partitioning 3-uniform hypergraphs},
Random Struct. Alg. {\bf 56} (2020) 1205--1221.

\bibitem{XYY2010}
B.~G. Xu, J. Yan and X. Yu,
\textit{Balanced judicious bipartitions of graphs},
J. Graph Theory {\bf 63} (2010) 210--225.

\bibitem{XY2011}
B.~G. Xu and X. Yu,
\textit{Better bounds for $k$-partitions of graphs},
Combin. Probab. Comput. {\bf 20} (2011) 631--640.

\bibitem{XY2014}
B.~G. Xu and X. Yu,
\textit{On judicious bisections of graphs},
J. Combin. Theory Ser. B {\bf 106} (2014) 30--69.

\bibitem{XYY2010}
B.~G. Xu, J. Yan and X. Yu,
\textit{Balanced judicious bipartitions of graphs},
J. Graph Theory {\bf 63} (2010) 210--225.

\bibitem{ZTY2015}
Y.~Zhang, Y.~C.~Tang and G.~Y.~Yan,
\textit{On judicious partitions of hypergraphs with edges of size at most 3},
European J. Combin. \textbf{49} (2015) 232--239.



\end{thebibliography}
\end{document}